\documentclass[12pt]{amsart}
\usepackage{amsmath,amscd,amssymb,amsfonts}
\newcommand{\h}{\hbox}
\newcommand{\q}{\quad}

\newcommand{\nin}{\par\noindent}
\newcommand{\bs}{\par\bigskip}
\newcommand{\ms}{\par\medskip}
\newcommand{\sk}{\par\smallskip}
\newcommand{\mopl}{\h{$\bigoplus$}}
\newcommand{\msum}{\h{$\sum$}}
\newcommand{\mcap}{\h{$\bigcap$}}
\newcommand{\al}{\alpha}
\newcommand{\C}{{\mathbf C}}
\newcommand{\D}{{\mathbf D}}
\newcommand{\De}{\Delta}
\newcommand{\Det}{\widetilde{\Delta}}
\newcommand{\Dt}{\widetilde{D}}
\newcommand{\Et}{\widetilde{E}}
\newcommand{\dd}{\partial}
\newcommand{\ft}{\tilde{f}}
\newcommand{\G}{{\mathcal G}}
\newcommand{\ga}{\gamma}
\newcommand{\gh}{\gamma_h}
\newcommand{\gc}{\gamma_c}
\newcommand{\Ht}{\widetilde{H}}
\newcommand{\jt}{{}\,\tilde{\!j}}
\newcommand{\K}{{\mathcal K}}
\newcommand{\LL}{{\mathbf L}}
\newcommand{\la}{\lambda}
\newcommand{\Omt}{\widetilde{\Omega}}
\newcommand{\ob}{{\mathbf 1}}
\newcommand{\PP}{{\mathbf P}}
\newcommand{\Q}{{\mathbf Q}}
\newcommand{\Sp}{{\rm Sp}}
\newcommand{\Ut}{\widetilde{U}}
\newcommand{\Xt}{\widetilde{X}}
\newcommand{\Yt}{\widetilde{Y}}
\newcommand{\Z}{{\mathbf Z}}
\newcommand{\Gr}{{\rm Gr}}
\newcommand{\Ker}{{\rm Ker}}
\newcommand{\into}{\hookrightarrow}
\newcommand{\bl}{\bigl}
\newcommand{\br}{\bigr}
\newcommand{\ssb}{\raise.15ex\h{${\scriptscriptstyle\bullet}$}}
\newcommand{\ssc}{\,\raise.15ex\h{${\scriptstyle\circ}$}\,}
\newcommand{\simto}{\buildrel{\sim}\over\longrightarrow}

\begin{document}
\title[Some remarks on limit mixed Hodge structure]
{Some remarks on limit mixed Hodge structure\\
and spectrum}
\author{Alexandru Dimca}
\address{Institut Universitaire de France et
Laboratoire J.A.\ Dieudonn\'e, UMR du CNRS 7351,
Universit\'e de Nice-Sophia Antipolis, Parc Valrose,
06108 Nice Cedex 02, France}
\email{Alexandru.DIMCA@unice.fr}
\author{Morihiko Saito}
\address{RIMS Kyoto University, Kyoto 606-8502 Japan}
\email{msaito@kurims.kyoto-u.ac.jp}
\maketitle
We give some remarks on limit mixed Hodge structure and spectrum.
These are more or less well-known to the specialists, and do not seem
to be stated explicitly in the literature. However, they do not seem
to be completely trivial to the beginners, and may be worth writing
down explicitly.
\bs\bs
\centerline{\bf 1.~Limit mixed Hodge structure}
\bs\nin
{\bf 1.1.} In [St1], [StZ], the limit mixed Hodge structures were
constructed in the unipotent monodromy case.
For the non-unipotent case, we can combine it with [St2] as follows.
Here we describe the limit of the mixed Hodge structure on the
cohomology with compact supports using the Cech-type construction,
since this seems to be the easiest way to explain the relation with
the theory of motivic Milnor fibers [DL].
For the usual cohomology (i.e.\ without compact supports), we can use
the commutativity of the dualizing functor $\D$ and the passage to the
limit mixed Hodge structure, i.e.\ $\D\ssc\psi_t=\psi_t\ssc\D$ (up to a
Tate twist), see [StZ], [Sai3].
Of course, we can also use the two weight filtrations on the
logarithmic complex associated with the divisor with $V$-normal
crossings [St2] as in [StZ]. It also follows from the theory of
mixed Hodge modules [Sai2], [Sai3].

Let $f:X\to\De$ be a projective morphism of complex manifolds where
$\De$ is an open disk.
We may assume that $f$ is smooth over $\De^*$ (shrinking $\De$ if
necessary).
Set $Y:=f^{-1}(0)$.
Let $D$ be a divisor on $X$ which is flat over $\De$, i.e.\ all the
irreducible components $D_j$ of $D$ are dominant over $\De$.
Assume $D\cup Y$ is a divisor with simple normal crossings. Set
$$U:=X\setminus D,\q f':=f|_U:U\to\De,\q D_J:=\mcap_{j\in J}D_j
\q(\h{where}\,\,\,D_{\emptyset}=X).$$
Let $Y_i$ be the irreducible components of
$Y\subset X$ with $m_i$ the multiplicity of $Y$ along the generic
point of $Y_i$.
Set $m={\rm LCM}(m_i)$.
Let $\ft:\Xt\to\Det$ be the normalization of the base change of
$f:X\to\De$ by the ramified $m$-fold covering $\pi_{\De}:\Det\to\De$
which is finite \'etale over $\De^*$, where $\Det$ is an open disk.
Let $\pi:\Xt\to X$ be the canonical morphism. Set
$$\Ut:=\pi^{-1}(U),\q\Yt:=\pi^{-1}(Y),\q\Dt:=\pi^{-1}(D),\q
\Dt_J:=\pi^{-1}(D_J).$$
Then $\Xt$ is a $V$-manifold, and $\Yt\cup\Dt$ is a divisor with
$V$-normal crossings on $\Xt$.
Let $\jt:\Ut\into\Xt$ be the natural inclusion.
There is a natural quasi-isomorphism
$$\jt_!\Q_{\Ut}\simto\K_{\Xt}^{\ssb}\q\h{with}\q
\K_{\Xt}^p:=\mopl_{|J|=p}\,\Q_{\Dt_J},
\leqno(1.1.1)$$
where the differential of $\K_{\Xt}^{\ssb}$ is defined in the same way
as a Cech complex as is well known.

Consider the complex
$$\psi_{\ft}\,\K_{\Xt}^{\ssb}.$$
This naturally underlies a cohomological mixed Hodge complex such that
its restriction to $\psi_{\ft}\,\Q_{\Dt_J}$ coincides with the one
defined in [St2] using the complex of logarithmic forms
$\Omt^{\ssb}_{\Dt_J}(\log(\Yt\cap\Dt_J))$ together with the Hodge
filtration $F$ and the weight filtration $W$ on it.
Indeed, we have canonical morphisms for $J\subset J'$
$$\Omt^{\ssb}_{\Dt_J}(\log(\Yt\cap\Dt_J))|_{\Dt_{J'}}\to
\Omt^{\ssb}_{\Dt_{J'}}(\log(\Yt\cap\Dt_{J'})),$$
and this is a filtered quasi-isomorphism for $W$ (forgetting the
filtration $F$), since $D\cup Y$ is a divisor with normal crossings.

There is a spectral sequence of mixed Hodge structures
$$_{\infty}E_1^{p,q}=\mopl_{|J|=p}\,H^q(\Dt_{J,\infty},\Q)
\Longrightarrow H^{p+q}_c(\Ut_{\infty},\Q),
\leqno(1.1.2)$$
which is induced by the truncations $\tau_{\ge k}$ on $\K_{\Xt}^{\ssb}$
for $k\in\Z$, and degenerates at $E_2$.
(This is the dual of the spectral sequence in [StZ], 5.7.)
Indeed, it is the `limit' by $\tt\to 0$ of the weight spectral sequence
$$_{\tt}E_1^{p,q}=\mopl_{|J|=p}\,H^q(\Dt_{J,\tt},\Q)\Longrightarrow
H^{p+q}_c(\Ut_{\tt},\Q),$$
where $\Xt_{\tt}:=\ft^{-1}(\tt)$ and $\Dt_{J,\tt}:=\Dt_J\cap\Xt_{\tt}$
for $\tt\in\Det^*$.
These spectral sequences are the dual of the spectral sequences in
[StZ] in the unipotent monodromy case.

\sk
Note that (1.1.2) is compatible with the actions of the semisimple
part $T_s$ and the nilpotent part $N:=(2\pi i)^{-1}\log T_u$ of the
monodromy $T$.

\ms\nin
{\bf 1.2.~The relation with motivic nearby fibers.}
With the above notation, let $E=Y\cup D$ with $E_i$ the irreducible
components of $E$.
We may assume $E_i=Y_i$ for $i\le r$ and $E_i=D_{i-r}$ for $i>r$,
where $r$ is the number of the irreducible components of $Y$.
For $I$ with $\min(I)\le r$ (i.e.\ $E_I\subset Y$), define
$$E_I=\mcap_{i\in I}E_i,\q
E_I^{\circ}=\mcap_{i\in I}E_i\setminus\mcap_{i\notin I}E_i,\q
\Et_I=\pi^{-1}(E_I),\q\Et_I^{\circ}=\pi^{-1}(E_I^{\circ}).
\leqno(1.2.1)$$
Note that $\Et_I^{\circ}\to E_I^{\circ}$ is a cyclic \'etale covering.
Set $I':=I\cap[1,r]$.
Let $\LL$ denote $\ob(-1)$ as an Chow motive where $\ob=[pt]$, and
$(-1)$ is the Tate twist, see e.g.\ [Mu], [Sch].
By [DL], [Lo] and [MT], [Ra], the motivic nearby fibers for the
morphisms $f:X\to\De$ and $f':U\to\De$ can be given respectively by
$$\msum_{\min(I)\le r}\,[(\Et_I^{\circ},T_s)](\ob-\LL)^{|I'|-1},\q
\msum_{\max(I)\le r}\,[(\Et_I^{\circ},T_s)](\ob-\LL)^{|I|-1}.
\leqno(1.2.2)$$
These belong to the Grothendieck group of Chow motives (with
$\Q$-coefficients) endowed with an action of $T_s$ of finite order by using equivariant resolutions of $(\Et_I,\Et_I\setminus\Et_I^{\circ})$,
see [DL].
Here $T_s$ denotes the semi-simple part of the monodromy, and is
given by the automorphism $\ga$ of $\Xt$ over $X$ induced by the base
change of the automorphism of $\Det$ defined by
$\tt\mapsto\zeta_m\,\tt$ with $\zeta_m:=\exp(2\pi i/m)$.
The action of $T_s$ on $\LL$ is the identity.
We denote the images of the two terms of (1.2.2) in the Grothendieck
group of mixed Hodge structures with an action of finite order
respectively by
$$\msum_{\min(I)\le r}\,[(H^{\ssb}_c(\Et_I^{\circ}),T_s]
(\ob-\LL)^{|I'|-1},\q
\msum_{\max(I)\le r}\,[(H^{\ssb}_c(\Et_I^{\circ}),T_s)]
(\ob-\LL)^{|I|-1},
\leqno(1.2.3)$$
where $H^{\ssb}_c(\Et_I^{\circ})$ is a complex of mixed Hodge structures
with zero differential, $\LL$ means here the class of $\Q(-1)$ with
trivial action of the monodromy, and $T_s$ is given by $(\ga^*)^{-1}$,
see (2.1.2) below.
Then these respectively coincide in the notation of (1.1.1) with
$$[(H^{\ssb}(\Xt_{\infty}),T_s)],\q[(H^{\ssb}_c(\Ut_{\infty}),T_s)].
\leqno(1.2.4)$$
Indeed, this follows from the construction of Steenbrink [St2]
together with the long exact sequence of mixed Hodge structures
$$\to H^j_c(Z')\to H^j_c(Z)\to H^j_c(Z\setminus Z')\to
H^{j+1}_c(Z')\to,
\leqno(1.2.5)$$
for any open immersions of complex algebraic varieties $Z'\into Z$,
which is compatible with the action of automorphisms of varieties.
(Here (1.2.5) can be proved by using mixed Hodge modules or the mapping
cone construction in [De3] together with the diagram of the octahedral
axiom of derived categories.)
The dual exact sequence of (1.2.5) for Borel-Moore homology is well
known in the theory of cycle maps of higher algebraic cycles.

Note that we get cohomology with compact supports in (1.2.4),
and this is quite different from the case of motivic Milnor fibers
in [DL].

\ms\nin
{\bf 1.3.~Remarks.} (i)
In case $E_I$ is simply connected, the cyclic \'etale covering
$\Et_I^{\circ}\to E_I^{\circ}$ can be determined by the multiplicities
$m_j$ of $Y$ along the irreducible components $Y_j$ intersecting $E_I$.
For example, assume $\max(I)\le r$ and
$$E_I=\mcap_{i\in I}Y_i=\PP^1,\q
E_I^{\circ}=E_I\setminus(E_{i'}\cup E_{i''})=\C^*,
\leqno(1.3.1)$$
with $i'\le r$ (i.e.\ $E_{i'}=Y_{i'}$). 
Then the covering degree of $\Et_I^{\circ}\to E_I^{\circ}$ and the
number of connected components of $\Et_I^{\circ}$ are given
respectively by
$${\rm GCD}(m_i\mid i\in I),\q{\rm GCD}(m_i\mid j\in I\cup\{i'\}).
\leqno(1.3.2)$$
This may simplify some argument in [MT].

\ms
(ii) In [DL], the semisimple part of the monodromy $T_s$ acts as an
automorphism of Chow motives.
This seems to be useful for the proof of the independence of the
motivic Milnor fiber by the resolutions of singularities.
For instance, we have $[\PP^1]=\ob\oplus\ob(-1)$  with
${\rm End}(\ob)={\rm End}(\ob(-1))=\Q$
in the category of Chow motives.
This is a special case of the Chow-K\"unneth decomposition,
see e.g.\ [Mu], [Sch].

\bs\bs
\centerline{\bf 2. Spectrum}
\bs\nin
{\bf 2.1.~Geometric monodromy and local system monodromy.}
Let $f:X\to S$ be a continuous map of topological spaces
which is locally topologically trivial over $S$.
We assume that the $H_j(X_s)$ and $H^j(X_s)$ with $\Q$-coefficients
are finite dimensional for any $j$.
Let $s\in S$, and $\ga\in\pi_1(S,s)$.
Let $\rho:Y\to[0,1]$ be the base change of $f$ by the loop $\ga$.
Choosing a trivialization over $[0,1]$, we get the
{\it geometric monodromy}
$$\ga_{\#}:X_s=Y_0\simeq Y_1=X_s,$$
where the middle homeomorphism is induced by the trivialization.
We have the induced action of the geometric monodromy on
homology and cohomology (with ${\bf Q}$-coefficients):
$$\ga_*\in{\rm Aut}(H_j(X_s)),\q\ga^*\in{\rm Aut}(H^j(X_s)),$$
such that
$$\ga^*={}^t\ga_*,
\leqno(2.1.1)$$
where $^t$ means the transpose.

On the other hand, we have the {\it local system monodromies}
$$\gh\in{\rm Aut}(H_j(X_s)),\q\gc\in{\rm Aut}(H^j(X_s)),$$
which are defined by using the following (trivial) local systems of
homology and cohomology groups over $[0,1]$:
$$\{H_j(Y_u)\}_{u\in[0,1]},\q\{H^j(Y_u)\}_{u\in[0,1]}.$$
The latter can be identified with the constant sheaf
$R^j\rho_*{\bf Q}_Y$, see also [De2], XIV, 1.1.2.
Note that $\gc$ coincides with the monodromy associated
to the nearby cycle functor if $f$ is a Milnor fibration.

The relations between the above monodromies are
given by
$$\ga_*=\gh,\q\ga^*=\gc^{-1}.
\leqno(2.1.2)$$
Indeed, the first assertion easily follows from the definition
(using simplicial chains for example).
We then get the second equality since
$$\ga^*={}^t\ga_*={}^t\gh=\gc^{-1},$$
where the last equality follows from
$$\langle \gc u,\gh v\rangle=\langle u,v\rangle\q\h{for}
\,\,\,u\in H^j(X_s),\,v\in H_j(X_s).
\leqno(2.1.3)$$
Here $\langle u,v\rangle$ denotes the canonical pairing between
cohomology and homology, and it can be extended to a canonical
pairing between the local systems so that (2.1.3) follows.

It does not seem that (2.1.2) has been clarified explicitly in the
literature.
In fact, it does not seem to cause big problems at least in the
local monodromy case since it is quasi-unipotent (except possibly for
the definition of spectrum as in [MT]).

\ms\nin
{\bf 2.2.~Example.}
Let $f$ be a homogeneous polynomial of $n$ variables with degree $d$,
having an isolated singularity at the origin.
Set $X=\C^n\setminus f^{-1}(0)$ and $S=\C^*$.
Here $f$ also denotes the morphism $X\to S$ induced by $f$.
Let $\ga$ be a generator of $\pi_1(S,s)\simeq{\bf Z}$
going around the origin counter-clockwise.
Then $\ga_{\#}$ is induced by the automorphism
$$\ga_{\#}:(x_1,\dots,x_n)\mapsto (\zeta_d\,x_1,\dots,\zeta_d\,x_n),
\leqno(2.2.1)$$
where $\zeta_d:=\exp(2\pi i/d)$, and $x_1,\dots,x_n$ are the
coordinates of $\C^n$.
The action of $\ga_{\#}$ is extended to an automorphism of
$\C[x_1,\dots,x_n]$ over $\C$ such that
$$\ga_{\#}^*x_i=\zeta_d\,x_i.
\leqno(2.2.2)$$
This can be checked for instance by
$\ga_{\#}^*(x_i-\zeta_d\,c_i)=\zeta_d(x_i-c_i)$.

\smallskip
Set $\omega=dx_1\wedge\cdots\wedge dx_n$, and
$$H''_f:=\Omega^n_{X,0}/df\wedge d\Omega^{n-2}_{X,0}.
\leqno(2.2.3)$$
This is called the Brieskorn lattice.
Let $g\in\C[x_1,\dots,x_n]$ be a monomial of degree $k$.
After Brieskorn, it is well known (and easy to show) that
$$\dd_tt(g\omega)=\h{$\frac{k+n}{d}$}g\omega\q\h{in}\,\,\,H''_f,
\leqno(2.2.4)$$
see e.g.\ the proof of Prop.~3.3 in [Sai1] for an argument in a
slightly more general case.

In the homogeneous polynomial case, we have moreover the well-known
relation
$$\gc=\exp(-2\pi i({\rm Res}\,t\dd_t)),
\leqno(2.2.5)$$
under the canonical isomorphism
$$H^{n-1}(X_1,\C)=H''_f/tH''_f,
\leqno(2.2.6)$$
where $s=1$.
The isomorphism can be defined by using a basis
$(\omega_1,\dots,\omega_{\mu})$ of $H''_f$ such that
$\dd_tt\omega_i=\al_i\omega_i$, and taking the restriction to $X_1$
after dividing the $\omega_i$ by $df$.
Indeed, the assertion is well known if the Brieskorn lattice is
replaced by the Deligne extension [De1].
In this case, the inverse isomorphism is given by
$$\h{$u\mapsto\exp(-\frac{\log t}{2\pi i}\log\gc)u,$}$$
for $u\in H^{n-1}(X_1)$ which is identified with a multivalued
section, and we have
$$\h{$t\dd_t\exp(-\frac{\log t}{2\pi i}\log\gc)u=
-\frac{\log\gc}{2\pi i}\exp(-\frac{\log t}{2\pi i}\log\gc)u,$}$$
where the eigenvalues of $-\frac{1}{2\pi i}\log\gc$ are chosen
corresponding to the Deligne extension.
This can be extended to the Brieskorn lattice case easily in the
homogeneous polynomial case.
So (2.2.5) follows.

By (2.2.4) and (2.2.5), the action of $\gc$ on $(g/df)|_{X_1}\in
H^{n-1}(X_1,\C)$ is given by the multiplication by
$$\exp(-2\pi i(k+n)/d).
\leqno(2.2.7)$$

On the other hand, (2.2.2) implies that the action of the geometric
monodromy (2.2.1) on $(g/df)|_{X_1}$ is given by the multiplication
by
$$\exp(2\pi i(k+n)/d).
\leqno(2.2.8)$$
This is the inverse of (2.2.7).

\ms\nin
{\bf 2.3.~Brieskorn lattices and mixed Hodge structures.}
The Brieskorn lattice $H''_f$ in (2.2.3) is defined for any
holomorphic function on a complex manifold $X$ having an isolated
singularity at $0\in f^{-1}(0)$.
It is a free $\C\{\!\{\dd_t^{-1}\}\!\}$-module of rank $\mu$, and is
contained in the Gauss-Manin system $\G_f$ which is the localization
of $H''_f$ by $\dd_t^{-1}$, i.e.\ $\G_f=H''_f[\dd_t]$.
The latter has the Hodge filtration defined by
$F^p\G_f:=\dd_t^{n-1-p}H''_f$ for $p\in\Z$, and also the filtration
$V$ of Kashiwara and Malgrange such that $\dd_tt-\al$ is nilpotent on
$\Gr_V^{\al}\G_f$.
By an argument similar to the proof of (2.2.6), there are isomorphisms
$$H^{n-1}(X_{f,0},\C)_{\la}=\Gr_V^{\al}\G_f\q\h{for}\,\,\,
\la=\exp(-2\pi i\al),
\leqno(2.3.1)$$
where $X_{f,0}$ is the Milnor fiber, and $V_\la$ denotes the
$\la$-eigenspace for any vector space $V$ with the action of the local
system monodromy $T$.
We have moreover
$$F^{n-1-q}H^{n-1}(X_{f,0},\C)_{\la}=\Gr_V^{\al}H''_f\q\h{for}\,\,\,
q<\al\le q+1,\,\,\,\la=\exp(-2\pi i\al),
\leqno(2.3.2)$$
where $F$ is the Hodge filtration of the mixed Hodge structure [St2],
see [ScSt], [Va], etc.
This is closely related with the definition of the spectrum in (2.4.5)
below.
In the case of Example~(2.2), it is related with [St3] by (2.2.4).

\ms\nin
{\bf 2.4.~Spectrum.}
Let $H$ be a mixed Hodge structure with a semisimple action $T$ of
finite order.
Set $H_{\C,\la}:=\Ker(T-\la)\subset H_{\C}$.
We define the spectrum $\Sp'(H,T)$ as in [Sai4] (and [DL]) by
$$\aligned&\Sp'(H,T):=\msum_{\al\in\Q}\,n'_{\al}t^{\al},\\
\h{with}\q &n'_{\al}=\dim_{\C}\Gr_F^pH_{\C,\la}\q
\h{for}\,\,\,p=[\al],\,\,\la=\exp(2\pi  i\al).\raise10pt\h{ }
\endaligned
\leqno(2.4.1)$$
For a holomorphic function $f$ on a complex manifold $X$ of dimension
$n$ and $x\in f^{-1}(0)$, we first define $\Sp'(f,x)$ by
$$\Sp'(f,x):=\msum_j\,(-1)^{n-1-j}\,\Sp'\bl(\Ht^j(X_{f,x}),T_s\br),
\leqno(2.4.2)$$
where $\Ht^j(X_{f,x})$ is the reduced Milnor cohomology endowed with
the canonical mixed Hodge structure, and $T_s$ is the semisimple part
of the local system monodromy $T$.
There are canonical isomorphisms
$$\Ht^j(X_{f,x})=H^ji_x^*\varphi_f\Q_X,
\leqno(2.4.3)$$
where $i_x:\{x\}\into X$ is the inclusion, see [De2].
They can be used to define the mixed Hodge structure on the left-hand
side.
Note that $T$ is equal to the inverse of the cohomological Milnor
monodromy by (2.1), and this is closely related with Example~(2.2) by
(2.3).

\sk
Let $\iota$ denote the involution of $\Z[t^{1/m},t^{-1/m}]$ over $\Z$
defined by
$$\iota(t^{\al})=t^{-\al}.$$
The spectrum $\Sp(f,x)$ is then defined by
$$\Sp(f,x):=t^n\iota(\Sp'(f,x)).
\leqno(2.4.4)$$
This spectrum $\Sp(f,x)$ coincides with the one in [St2] for the
isolated singularity case (using the complex conjugate of (2.4.1)
together with the symmetry (2.4.6) below).
It coincides with the one in [St4] up to the
multiplication by $t$.
Indeed, the above definition of $\Sp(f,x)$ can be rewritten as
$\Sp(f,x):=\msum_{\al\in\Q}\,n_{\al}t^{\al}$ with
$$\aligned&n_{\al}=\msum_j\,(-1)^{n-1-j}\,\dim_{\C}\Gr_F^{n-1-q}
\Ht^j(X_{f,x},\C)_{\la}\\
&\h{for}\,\,\,q<\al\le q+1,\,\la=\exp(-2\pi  i\al),\endaligned
\leqno(2.4.5)$$
and this is used in {\it loc.~cit.} (up to the multiplication by $t$).
In the isolated singularity case, the formula (2.4.5) is closely
related with (2.3.2) and also with the calculations in (2.2).

If $f$ has an isolated singularity at $x$, we have the symmetry of
mixed Hodge numbers by [St2] so that
$$\Sp(f,x)=\Sp'(f,x).
\leqno(2.4.6)$$

In case $f$ is a weighted-homogeneous polynomial of weights
$(w_1,\dots,w_n)$ (i.e.\ $f$ is a linear combination of monomials
$x_1^{m_1}\cdots x_n^{m_n}$ with $\sum_i w_im_i=1$) and has an
isolated singularity at $0$, it is well known that
$$\Sp(f,0)=\prod_{i=1}^n\biggl(\frac{t-t^{w_i}}{t^{w_i}-1}\biggr).
\leqno(2.4.7)$$
This follows from [St3], see also [St2] and the proof of
Proposition~5.2 in [Di].
In the case of homogeneous polynomials (i.e.\ $w_i=1/d$ for any $i$),
this follows also from the calculation in Example~(2.2) using (2.3).

For a polynomial mapping $f:\C^n\to\C$, we can define the spectrum at
infinity (see also [Sab1]) by
$$\aligned\Sp'(f,\infty)&:=\msum_j\,(-1)^{n-1-j}\,
\Sp'(\Ht^j(X_{\infty}),T_s),\\
\Sp(f,\infty)&:=t^n\iota(\Sp'(f,\infty)),\endaligned
\leqno(2.4.8)$$
where $\Ht^j(X_{\infty})$ is the limit mixed Hodge structure of
$\Ht^j(X_t)$ at infinity (of $\C$), and $T_s$ is the semisimple part
of the local system monodromy $T$ associated with a sufficiently large
loop around the origin which goes counter-clockwise from the origin
(and clockwise from $\infty\in\PP^1$).
This definition is compatible with the one in the weighted-homogeneous
isolated singularity case in (2.4.7).
In the cohomologically tame case, we have the symmetry by [Sab2]
(i.e.\ Theorem~1 in this paper) so that (2.4.6) holds, and the
definition~(2.4.8) seems to coincide with the one in [MT] (where
the cohomology with compact supports is used) if the local system
monodromy is used there.


\begin{thebibliography}{Sab2}

\bibitem[De1]{De1}
Deligne, P., Equation diff\'erentielle \`a points singuliers
r\'eguliers, Lect.\ Notes in Math.\ 163, Springer, Berlin, 1970.

\bibitem[De2]{De2}
Deligne, P., Le formalisme des cycles \'evanescents, in SGA7 XIII
and XIV, Lect. Notes in Math. 340, Springer, Berlin, 1973,
pp. 82--115 and 116--164.


\bibitem[De3]{De3}
Deligne, P., Th\'eorie de Hodge, III, Publ.\ Math.\ I.H.E.S. 44
(1974), 5--77.

\bibitem[DL]{DL}
Denef, J.\ and Loeser, F., Motivic Igusa zeta functions,
J.\ Alg.\ Geom.\ 7 (1998), 505--537.

\bibitem[Di]{Di}
Dimca, A., Monodromy and Hodge theory of regular functions,
in New Developments in Singularity Theory, Kluwer Acad.\ Publ.,
Dordrecht, 2001, pp. 257--278.

\bibitem[Lo]{Lo}
Looijenga, E., Motivic measures, in S\'eminaire Bourbaki 1999/2000,
Ast\'erisque 276 (2002), 267--297.

\bibitem[MT]{MT}
Matsui, Y.\ and	Takeuchi, K.,
Monodromy at infinity of polynomial maps and Newton polyhedra,
preprint (arXiv:0912.5144v11).

\bibitem[Mu]{Mu}
Murre, J.P., On the motive of an algebraic surface, J. reine angew.\
Math.\ 409 (1990), 190--204.

\bibitem[Ra]{Ra}
Raibaut, M., Fibre de Milnor motivique \`a l'infini,
C.\ R.\ Acad.\ Sci.\ Paris, Ser.\ I 348 (2010), 419--422.

\bibitem[Sab1]{Sab1}
Sabbah, C., Monodromy at infinity and Fourier transform,
Publ. RIMS, Kyoto Univ.\ 33 (1997), 643--685.

\bibitem[Sab2]{Sab2}
Sabbah, C., Hypergeometric periods for a tame polynomial,
Port.\ Math., 63 (2006), 173--226 (or arXiv:math/9805077).

\bibitem[Sai1]{Sai1}
Saito, M., Exponents and Newton polyhedra of isolated hypersurface
singularities, Math.\ Ann.\ 281 (1988), 411--417.

\bibitem[Sai2]{Sai2}
Saito, M., Modules de Hodge polarisables, Publ. RIMS, Kyoto Univ.\
24 (1988), 849--995.

\bibitem[Sai3]{Sai3}
Saito, M., Mixed Hodge modules, Publ. RIMS, Kyoto Univ.\ 26
(1990), 221--333.

\bibitem[Sai4]{Sai4}
Saito, M., On Steenbrink's conjecture, Math.\ Ann.\ 289 (1991),
703--716.

\bibitem[ScSt]{ScSt}
Scherk, J.\ and Steenbrink, J.H.M., On the mixed Hodge structure on
the cohomology of the Milnor fibre, Math.\ Ann.\ 271 (1985), 641--665.

\bibitem[Sch]{Sch}
Scholl, A.J. , Classical Motives, Proc.\ Symp.\ Pure Math.\ 55
(1994), Part 1, 163--187.

\bibitem[St1]{St1}
Steenbrink, J.H.M., Limits of Hodge structures, Inv.\ Math.\ 31
(1976), 229--257.

\bibitem[St2]{St2}
Steenbrink, J.H.M., Mixed Hodge structure on the vanishing cohomology,
in Real and complex singularities, Sijthoff and Noordhoff, Alphen aan
den Rijn, 1977, pp. 525--563.

\bibitem[St3]{St3}
Steenbrink, J.H.M., Intersection form for quasi-homogeneous
singularities, Compositio Math.\ 34 (1977), 211--223.

\bibitem[St4]{St4}
Steenbrink, J.H.M., The spectrum of hypersurface singularities,
Ast\'erisque 179-180 (1989), 163--184.

\bibitem[StZ]{StZ}
Steenbrink, J.H.M.\ and Zucker, S., Variation of mixed Hodge structure
I, Inv.\ Math., 80 (1985), 489--542.

\bibitem[Va]{Va}
Varchenko, A.N., Asymptotic Hodge structure in the vanishing
cohomology, Math.\ USSR Izv.\ 18 (1982), 469--512.

\end{thebibliography}
\end{document}